\documentclass{article}
\usepackage[mathcal,mathbf]{euler}
\usepackage{inputenc, theorem,amsmath,enumerate,fancyhdr,amssymb,amsfonts,accents, url, color, float,graphicx, subfigure, hyperref}
\usepackage{algorithm, algpseudocode}

\algnewcommand{\MyStateA}[1]{\State \parbox[t]{\dimexpr\linewidth-\algorithmicindent}{\hangindent=\algorithmicindent \hangafter=1 \strut #1\strut}}
\algnewcommand{\MyStateB}[1]{\State \parbox[t]{\dimexpr\linewidth-\algorithmicindent-\algorithmicindent}{\hangindent=\algorithmicindent \hangafter=1 \strut #1\strut}}
\usepackage{MyCommands}

\title{A branch-and-bound algorithm for the minimum radius $k-$enclosing ball problem}

\author{Marta Cavaleiro\footnote{Rutgers University, MSIS Department \& RUTCOR, 100 Rockafellar Rd, Piscataway, NJ 08854. \texttt{marta.cavaleiro@rutgers.edu} } \and Farid Alizadeh\footnote{Rutgers University, MSIS Department \& RUTCOR, 100 Rockafellar Rd, Piscataway, NJ 08854. \texttt{farid.alizadeh@rutgers.edu} } }

\date{\today} 

\begin{document}

\maketitle

\abstract{The minimum $k$-enclosing ball problem seeks the ball with smallest radius that contains at least~$k$ of~$m$ given points in a general $n$-dimensional Euclidean space. This problem is NP-hard. We present a branch-and-bound algorithm on the tree of the subsets of~$k$ points to solve this problem. The nodes on the tree are ordered in a suitable way, which, complemented with a last-in-first-out search strategy, allows for only a small fraction of nodes to be explored. Additionally, an efficient dual algorithm to solve the subproblems at each node is employed.}

\bigskip


\section{Introduction} 
For a set $\Pp:=\{p_1,\hdots,p_m\}$ of $m$ points in a general $n$-dimensional Euclidean space, define a $k$-enclosing ball of $\Pp$ as an Euclidean hypersphere that covers at least $k$ points of $\Pp$. The \emph{minimum $k$-enclosing ball} problem,  $MkEB$ in short, seeks the $k$-enclosing ball of $\Pp$ with smallest radius. It can be formulated as
\begin{equation}\label{mkeb}
MkEB(\Pp):=\begin{array}{cl}\displaystyle
\min_{r, x} & r  \\
\text{s.t.} & \left\lvert\left\{p_j\in \mathcal P:\, \|x-p_j \| \leq r \right\} \right\rvert\geq k,
\end{array}
\end{equation}
\noindent with $|\,.\,|$ denoting the cardinality of a set and $\|\,.\,\|$ denoting the Euclidean norm. When $k=m$, the problem becomes the well known \emph{minimum enclosing ball} (MEB) problem, also called the $1$-center problem:
\begin{equation} \label{meb}
MEB(\Pp) : = \begin{array}{cl}\displaystyle
\min_{r, x} & r \\
\text{s.t.} & \|x-p_j\|\leq r, \quad p_j \in \Pp.
\end{array}\end{equation}

\noindent The MkEB problem can then be seen as the robust version of the MEB problem. In many applications, data is usually noisy, and the presence of outliers can affect the solution of the MEB problem dramatically, justifying the need for a robust solution. Both are fundamental problems in computational geometry having many applications in areas like data mining, statistics, image processing, machine learning, etc. 

\smallskip

 The initial approaches to solve the MkEB problem consisted on constructing higher order Voronoi diagrams and then performing a search in all or some of the Voronoi cells (see Aggarwal et al. in~\cite{Aggarwal91}). A variety of methods were developed next: Efrat et al. \cite{Efrat93} developed an algorithm to solve the problem using the parametric search technique (due to Megiddo \cite{Megiddo83}); Eppstein and Erickson~\cite{Eppstein93} started by computing the $\mathcal{O}(k)$ nearest neighbors for each point, reducing the initial problem to $\mathcal{O}(m/k)$ smaller subproblems (see also \cite{Datta95} for an improvement on their algorithm); and Matou\v{s}ek developed a simple randomized search algorithm in \cite{Matousek95}. These methods solved the problem exactly, and though they were targeting the planar case, their extension to higher dimensions is possible.

Approximation algorithms started to be developed more recently. Har-Peled and Mazumdar \cite{HarPeled05} focused on the planar case and presented a linear $2$-approximation algorithm by dividing the plane in \emph{grids}, based on which they developed both a randomized exact algorithm and a $\epsilon$-approximation algorithm for the MkEB problem. The same ideas were further explored in \cite{HarPeled15} where an $\epsilon$-algorithm with expected $\bigO(m/\epsilon^n)$ time was shown as an application of a more general framework to solve several computational geometry problems. Other approximation algorithms have been developed based on the computation of \emph{(robust) core sets} \cite{HarPeled04, Agarwal06}.

Recently, Shenmaier \cite{Shenmaier13} proved that the minimum $k$-enclosing ball problem is NP hard in the strong sense, by reducing in polynomial time the NP-complete $k$-clique problem in a regular graph to the MkEB problem. In the same paper, a polynomial-time $\epsilon$-approximation algorithm was also proposed with complexity $\bigO(m^{1+1/\epsilon^2}n)$.

\smallskip

The minimum enclosing ball problem (\ref{meb}) can be easily formulated as a second-order-cone program (SOCP), more specifically, as a quadratic program (QP) \cite{Gartner00}, and thus solved in polynomial time. For a convex quadratic program, that is a problem with a convex quadratic objective and polyhedral constraints, simplex-type methods have been known for decades, see for instance \cite{Wolfe59, Panne64, Goldfarb83}. However, due to its particular geometric structure, many special purpose simplex-type methods have been developed to solve the MEB problem. Among them there are both primal and dual iterative algorithms that have an analogous mechanics to the simplex method for linear programming \cite{Fischer03, Dearing09}.

\smallskip

In this paper, we present  a branch-and-bound (B$\&$B) algorithm that builds the tree of all $k$-subsets of $\Pp$, that is, subsets of size $k$, and that uses the algorithm presented in \cite{Cavaleiro17}, an improvement on \cite{Dearing09}, to solve the minimum enclosing ball problem in each node. Branch-and-bound methods using the same search tree have also been considered in \cite{Candela96} to solve the minimum $k$-enclosing ellipsoid, and also in \cite{Narendra77, Yu93, Chen03} to address the feature selection problem.

\smallskip

The paper is organized as follows. In section \ref{sec:BBdescription}, we introduce the basics of the B$\&$B algorithm. Then, in section \ref{sec:BBfeatures}, we describe how the B$\&$B can be modified to further improve the performance and provide a detailed description of the revised algorithm. In section \ref{sec:MEB} we cover the algorithm chosen to solve the subproblem at each node of the tree. In section \ref{sec:experiments}, we give experimental results, and finally, draw conclusions in section \ref{sec:conclusions}.

\section{The branch-and-bound algorithm}\label{sec:BBdescription}
\begin{figure}[H]
	\centering\includegraphics[scale=0.25]{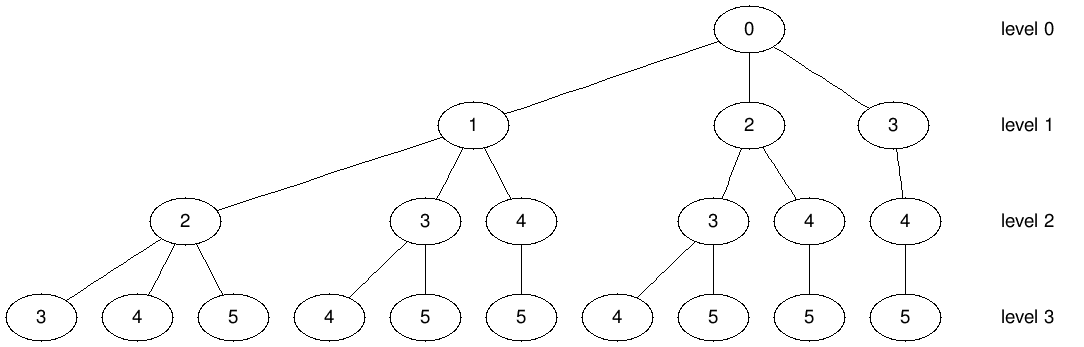}
	\caption{Search tree for $m=5$ and $k=3$.}
	\label{fig_tree}
\end{figure}

Consider the minimum $k$-enclosing ball problem (\ref{mkeb}). The search space, that is the set of all $k$-subsets of~$\Pp$, can be represented as a tree, as shown in Figure \ref{fig_tree} for the case of $m=5$ and $k=3$.

\noindent Each node on the tree corresponds to one point of $\Pp$ (except the root, which represents the empty set). The set of nodes along any path from the root to a certain node represents a different subset of~${\Pp}$ with as many points as the tree level number. Thus, each path to a leaf represents a different $k$-subset of~${\Pp}$. 

\smallskip

The branch-and-bound  algorithm constructs the tree in a top-down manner, as it searches for the minimum radius $k$-enclosing ball of $\Pp$. At each node $N$ of the tree, the algorithm solves the following subproblem: find the minimum enclosing ball that encloses the points corresponding to the nodes on the path from the root to $N$. The algorithm does not need necessarily to reach a leaf to find a ball that encloses at least $k$ points, since, at any node, there may be other points of ${\Pp}$ that are enclosed besides those on the respective path. Information about the radius, $r^*$, and center, $x^*$, of the currently smallest $k$-enclosing ball found so far is kept at all times. We call $r^*$ the \emph{bound}. 

The search begins at level $0$, when the root node is created and added to the pool of \emph{live} nodes~$\Ll$. A \emph{live} node is a node that has been explored, that is, whose corresponding MEB problem has been solved, but whose immediate successors (child nodes) have not been generated yet. At each iteration, and following some selection rule, a node is removed from ${\Ll}$ to be branched. Let $r$ and $x$ be the solution to the node's corresponding MEB problem. If, for that node, $r\geq  r^*$, then its successor nodes do not need to be analyzed, since the radius can never decrease as we proceed to higher levels down the tree, since a new point is being added at each level. On the other hand, if $r< r^*$, each child node is created and the associated MEB problem solved. If the radius $r$ corresponding to a child node is such that $r\geq  r^*$, then that node can be immediately disregarded. Otherwise, if the child's MEB encloses~$k$ points of~$\Pp$, then its radius becomes the new bound. If neither of these situations occur, the child node is simply added to $\Ll$ to be branched later. The algorithm terminates when $\Ll$ is empty, meaning that all nodes were either analyzed or cut off from the tree. The solution of (\ref{mkeb}) is the radius $r^*$ and center $x^*$ as defined at termination.

\section{Features of the B\&B}\label{sec:BBfeatures}
\subsection{Search tree design}\label{subsec:treeDesign}

It is evident that, given a search tree topology, different point-node assignments may be defined. Our B$\&$B algorithm aims at positioning ``bad''  points, points that will likely give a higher bound, on the denser part of the tree (left side), and, ``good'' points on the less dense part of the tree (right side). Nodes on the denser part are therefore more likely to be cut off.

We introduce a node-ordering method that works as follows. When a node $N$ is selected for branching, first, the distances between the center of its MEB and each of the points corresponding to points on node $N$'s subtree are calculated and sorted in descending order. The $j$-th child node of $N$ corresponds to the point with the $j$-th largest distance. This in-level node-ordering heuristic is effective because, for a general set $\Q$ and a point $q$, the distance from $q$ to the center of $MEB(\Q)$ is a good predictor of the radius of the $MEB(\Q\cup \{q\})$. Following this heuristic, nodes with higher radius will tend to be on the denser side of the tree, resulting in a better chance of larger-scale pruning. As expected, when compared to simply ordering the nodes in lexicographic order, or randomly, this heuristic caused the number of evaluated nodes to decrease considerably.

\smallskip

One could also think of picking a node's children as explained before, and sorting them, after solving the MEB problem for each, by the corresponding radius. While this approach could result in a (small) decrease on the number of explored nodes, such benefit is outweighed by the increase in the running time due to an additional sorting at each node.

\smallskip

\emph{Children of the root node:} The $m-k+1$ child nodes of the root are special and the heuristic above cannot be applied to them. An approach that proved to give good results in practice consists of first calculating the minimum enclosing ball of all points in $\Pp$, and then selecting the $m-k+1$ farthest points from its center as the children of the root. 

\smallskip

\begin{figure}[H]
	\centering \includegraphics[scale=0.25]{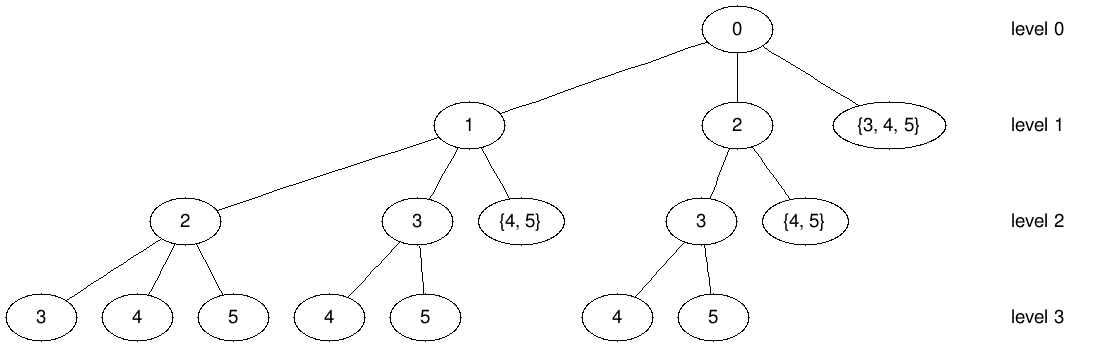}
	\caption{Minimum solution search tree for $m=5$ and $k=3$.}
	\label{fig_mintree}
\end{figure}

\emph{The minimum solution tree\footnote{The concept of minimum solution tree was first introduced by Yu and Yuan in \cite{Yu93} for a branch-and-bound algorithm for solving feature selection problems.}:} Consider Figure \ref{fig_tree}. The rightmost nodes of any subtree are all degree-one nodes and therefore simply constitute a path. Instead of evaluating each node on such paths one at a time, we can jump directly to the leaf and solve the corresponding MEB problem. The path below such one-degree nodes can therefore be merged into a \emph{special} single node, yielding a modified topology of the search tree called the \emph{minimum solution tree}, as shown in Figure~\ref{fig_mintree}. 

\noindent The computational work of solving the MEB problem at a \emph{special} node can only be less than the one of solving separately the MEB for each one of the nodes it contains. However, this improvement is not expected to be significant due to the specifics of the algorithm chosen to solve the MEB problem at each node (see Section \ref{sec:MEB}). Nevertheless, the use of the minimum solution tree leads to time savings in our implementation since it allows to avoid the addition and removal from $\Ll$ of nodes from the aforementioned paths.

\subsection{Search strategy}\label{subsec:search}

Assume that, when branching on a node, its child nodes are created and added to $\Ll$ in a left-to-right order. $\Ll$ is maintained as a stack, that is, the last node added is the first one to be removed. This last-in-first-out (LIFO) strategy offers several advantages. First, a stack is an easy data structure to maintain. Second, given our node-ordering heuristics, this LIFO strategy implies that the less dense part of the tree, that is the most promising part, is searched first. Moreover, unless a leaf is reached or a node is cut off (because it leads to a radius larger than the current bound), at every iteration the B$\&$B descends one level. This underlying depth-first-search allows reaching leaves fast. Finally, it is possible to calculate the maximum size of $\Ll$, as Theorem \ref{th_maxL} states, allowing the required memory to be allocated before hand.

\begin{theorem} [Maximum length of $\Ll$]\label{th_maxL}
	When LIFO is used as node-selection rule and the minimum solution tree is adopted, the maximum length of $\Ll$ is $m-k$. This upper bound is tight.
\end{theorem}

\proof{ Let $l(N)$ be the level of node $N$ and $N_0, N_1,\hdots ,N_{l(N)}$ be the nodes on the path from the root $N_0$ to $N\equiv N_{l(N)}$. Let us define the index of~$N$, $i(N)$, as the index of the point that~$N$ would correspond to if the nodes were to be organized in lexicographical order. By definition~$i(N_0)=0$. Together, $i(N)$ and $l(N)$ uniquely define the shape of $N$'s subtree. It is easy to see that the number of immediate successors of any node $N$ is $m-k+l(N)-i(N)+1$ and that $i(N)\leq m-k+l(N)$ always holds.
	
Suppose the worst case scenario: no subtree is ever cut-off, that is, every node will be added to $\Ll$ at some point, with the exception of the leaves (which are disregarded after being evaluated). 

Consider any stage of the algorithm when a node $N$ was removed from $\Ll$ to be branched. Note that $N$ cannot be a leaf, therefore $l(N)\leq k-1$. If $l(N)=0$, then $N$ is the root and $\Ll$ is empty. The root has $m-k+1$ immediate successors, each of which will be added to $\Ll$ with the exception of the last one, which is a leaf since the minimum solution tree is being adopted. After this process $\Ll$ has size $m-k$. 

Consider now the case $l(N) \geq 1$.  Since $N$ was removed from $\Ll$, all nodes on the path from the root to $N$ have also been removed from $\Ll$. Moreover, $\Ll$ cannot contain any nodes of higher level than~$l(N)$. Let $j$ be any level between $1$ and $l(N)$. The indexes of the nodes of level $j$ that are still in $\Ll$ are $i(N_{j-1})+1, \hdots, i(N_j)-1$. Thus, the number of nodes in $\Ll$ is:

\begin{equation*}
\sum_{j=1}^{l(N)}\left(i(N_j)- i(N_{j-1})-1\right)= i(N)-l(N).
\end{equation*}

\noindent Once $N$ is branched, its $m-k+l(N)-i(N)+1$ immediate successor nodes are created and evaluated, and with the exception of the last one, they are all added to $\Ll$. The size of $\Ll$ then becomes $m-k$.

}

\medskip

We experimented with other search strategies but the LIFO showed to be the most efficient. Breadth-first-search and best-first-search strategies tend to keep a very large number of nodes in $\Ll$, which leads to substantial memory space usage, and, in the case of best-first-search, more computational effort maintaining the priority queue. Note that, as a consequence of our node-ordering method, the LIFO strategy can be seen as a combination of depth-first-search and best-first-search.

\subsection{Getting an initial solution}\label{subsec:initSol}
A good initial solution, a ball that covers at least $k$ points and with a radius close to the optimal one, is expected to result in more aggressive pruning. There are many approaches to get a good initial $k$-enclosing ball. For instance, we could randomly select $k$ points of~$\Pp$ and find their minimum radius enclosing ball. Or, we could pick a random point in~$\Pp$ and find the~$k-1$ closest points to it, and then find the minimum radius ball of the resulting set. However, in general, these methods do not produce a ball with an enough small radius. Other, usually more successful, approaches are described next. They were also studied in \cite{Ahipa15} for the more general problem of the minimum volume ellipsoid with partial inclusion: 

\emph{Spherical ordering:} First, find the minimum enclosing ball of $\Pp$. Next, pick the $k$ closest points to the center of the ball. Finally, find the minimum enclosing ball of those $k$ points.

\emph{Spherical peeling:} Find the minimum enclosing ball of $\Pp$, discard one of the points from the boundary, and find the minimum enclosing ball of the remaining points. Repeat this until only $k$ points remain (note that the final ball may still cover more than $k$ points).

\smallskip

If the points have a ``spherical'' shape 
and either their density is uniform (e.g. points sampled uniformly from a ball) or there is a high concentration of points closer to the center of the sphere (e.g. points from a standard normal distribution), then spherical ordering naturally gives a better approximation for the minimum $k$-enclosing ball independently of the value of~$k$. On the other hand, when the density of points is low near the center, but the data is still ``spherical'' (e.g. points sampled uniformly from the surface of a sphere) the best approach depends on $k$. When $k$ is large then spherical ordering performs well, but when $k$ is small, spherical ordering will likely return a ball centered close to the center of the data and that covers a very sparse set of points, resulting in a bad initial solution. In this situation, randomly selecting a point and finding the minimum enclosing ball of its $k-1$ closest points is likely to provide a better initial solution.

When the data does not have a spherical shape, spherical peeling is likely to give better approximations, since, by iteratively removing points from the boundary of the minimum enclosing ball of the remaining data, the final ball will likely cover a dense subset of $k$ points.

The choice of the method to obtain an initial solution should therefore be based on the arrangement of the points in $\Pp$ and on $k$, whenever possible.

\smallskip

The spherical ordering and spherical peeling methods may instead be implemented by starting with the minimum enclosing ball of any subset of $\Pp$ with $k$ points (picked randomly, for instance). In practice, we observed that this random approach provides just as good initial solutions as the methods which start with all of $\Pp$.

\smallskip

In Section \ref{subsec:5is} we present some studies on the maximum improvement a good initial solution can provide for different datasets and values of $k$.

\subsection{Lower bounds}\label{subsec:LB}

Obtaining a {lower bound} for the radius of the MkEB  is useful for cases where we are satisfied with suboptimal solutions, or cases where there is a limit on the number of iterations, or the running time  of the algorithm. Having a good lower bound will indicate approximately how good our current best solution is.

A lower bound can be easily obtained by observing that the diameter of the MEB of any set of points $\Q$ (twice the radius of the ball) is larger than or equal to $\diam(\Q)$, the diameter of set $\Q$ (the maximum Euclidean distance between any two points in $\Q$). Suppose the solution of (\ref{mkeb}) is a ball which encloses the subset $\Q \subseteq \Pp$ and has radius $r^*$ (notice that $|\Q|\geq k$). We then have

\begin{equation}\displaystyle\begin{array}{ccl}\label{lbdist}
	2r^* &\geq& \diam(\Q)\\
	&=&\left(\frac{|\Q|(|\Q|-1)}2\right)\text{-th smallest pairwise distance of $\Q$}\\
	& \geq& \left(\frac{|\Q|(|\Q|-1)}2\right)\text{-th smallest pairwise distance of $\Pp$}\\ 
	&\geq& \left(\frac{k(k-1)}2\right)\text{-th smallest pairwise distance of $\Pp$}.
\end{array}\end{equation}

\noindent A lower bound for $r^*$ is therefore half the $\left(k(k-1)/2\right)$-smallest pairwise distance between any two points of~$\Pp$.

An alternative way to obtain a lower bound for the optimal radius is to solve a relaxation of~(\ref{mkeb}). It is possible to formulate~(\ref{mkeb}) as a mixed-integer-second-order-cone problem as follows:

\begin{equation}\begin{array}{cl}\label{lbSOCP}
			\displaystyle\min_{r,x,\beta_j} & r  \\
		\text{s.t.} & \left\|\beta_jp_j+(1-\beta_j)v_j-x\right\|\leq r, \quad p_j\in \Pp\\
		& \displaystyle\sum_{j=1}^{m} \beta_j \geq k\\
		&\beta_j\in \{0,1\},\, j=1,...,m,
	\end{array}\end{equation}
	
\noindent where $v_j$ is any point that is guaranteed to be covered by the optimal $k$-enclosing ball of $\Pp$. If we relax the integrality constraints of (\ref{lbSOCP}) to $0\leq\beta_j\leq 1$, and solve the resulting SOCP $m-k+1$ times, each time with $v_j$ being a different point of $\Pp$, we obtain a lower bound for $r^*$ by selecting the minimum radius among all the obtained solutions. The relaxation of (\ref{lbSOCP}) is not known to be equivalent to a quadratic programming problem, so as a result we are not aware of a simplex-type algorithm to solve it. Instead, interior point methods seem to be the only algorithms that can solve it effectively. Consequently, setting up a branch-and-bound algorithm centered on this relaxation does not seem to be efficient, since using interior point methods at each node is excessively time consuming.
	
Problem (\ref{mkeb}) can also be formulated as a mixed-integer-quadratic problem as follows:

\begin{equation}\begin{array}{cl}\label{lbQP}
		\displaystyle	\min_{r,x,\beta_j} & r^2  \\
		\text{s.t.} & \|p_j-x\|^2\leq r^2+(1-\beta_j)M, \quad p_j\in \Pp\\
		& \displaystyle\sum_{j=1}^{m} \beta_j \geq k\\
		&\beta_j\in \{0,1\},\, j=1,...,m.
	\end{array}\end{equation}
		
\noindent Relaxing the integrality constraints yields a lower bound for $r^*$. However, the value of $M$ needs to be very large so as to not jeopardize feasibility. Without any knowledge about the location of the minimum $k$-enclosing ball in $\Pp$, all we can conclude is that $M$ is smaller than the square of~$\diam(\Pp)$, which could be quite large. As a result, the relaxation ($\ref{lbQP}$) tends not to be very useful, except for when $k$ is very close to the number of points in $\Pp$.
	
\medskip

The lower bounds for the solution of problem (\ref{mkeb}) presented above can be modified to calculate, at each node $N$ of the search tree, a lower bound for the best solution that can be found on the nodes of its subtree. When that lower bound is higher than the smallest radius of a $k$-enclosing ball found so far, node $N$ can be cut off. Consider node $N$ of the search tree. Let $\Ss t(N)$ be the points of $\Pp$ corresponding to the nodes on $N$'s subtree, and $\Pp(N)$ be the points of $\Pp$ corresponding to nodes on the path from the root node to $N$ (included). The smallest radius of a $k$-enclosing ball that can be found in any of $N$'s subtree nodes is given by the problem of finding the minimum radius $k$-enclosing ball that covers $\Pp(N)$ and at least $k-l(N)$ points of $\Ss t(N)$. That problem can be formulated based on (\ref{lbSOCP}):

\begin{equation}\begin{array}{cll}\label{lbSOCPnode}
			\displaystyle\min_{r,x,\beta_j} & r  &\\
		\text{s.t.} & \|p_i-x\|\leq r, & p_i\in \Pp(N)\\
		& \|\beta_jq_j+(1-\beta_j)v_j-x\|\leq r, & q_j\in \Ss t(N)\color{black}\\
		& \displaystyle\sum_{j=1}^{|\Ss t(N)|} \beta_j \geq k-l(N)&\\
		&\beta_j\in \{0,1\},& j=1,...,|\Ss t(N)|.
\end{array}\end{equation}
	
\noindent where $v_j$ is a point that is guaranteed to be covered by the smallest $k$-enclosing ball found on $N$'s subtree. Relaxing the integrality condition on the multipliers $\beta_j$ yields a SOCP that provides a lower bound. Picking $v_j$ such that the distance between $v_j$ and $q_j$ is small gives the best bound, although finding such $v_j$ causes an increase in the running time. A better alternative is choosing~$v_j$ as the closest point of $\Pp(N)$ to~$q_j$.
	
A formulation based on (\ref{lbQP}) is obtained in a similar way:
\begin{equation}\begin{array}{cll}\label{lbQPnode}
		\displaystyle	\min_{r,x,\beta_j} & r^2 & \\
			\text{s.t.} & \|p_i-x\|^2\leq r^2, & p_i\in \Pp(N)\\
			& \|q_j-x\|^2\leq r^2+(1-\beta_j)M, & q_j\in \Ss t(N)\\
			& \displaystyle\sum_{j=1}^{| \Ss t(N)|} \beta_j \geq k-l(N)&\\
			&\beta_j\in \{0,1\},& j=1,...,|\Ss t(N)|.
\end{array}\end{equation}
		
\noindent As usual, the relaxation of (\ref{lbQPnode}) can be solved to obtain a lower bound.

A lower bound on the best solution of a node's subtree can be obtained based on (\ref{lbdist}). However, this lower bound is not advantageous given the node-ordering method we adopt. The points of $\Pp(N)$ are among the ones with higher pairwise distances of the set $\Pp(N)\cup \Ss t(N)$, given the way they were selected. Therefore, in most cases, half the $(k(k-1)/2)$-smallest distance between the points in $\Pp(N)\cup \Ss t(N)$ will be smaller than the radius of the MEB of $\Pp(N)$.

\smallskip

Finally, we point out that these lower bounds on the best solution of a node's subtree are only useful when the number of points in $\Ss t(N)$ is large, so it is not efficient to calculate such lower bounds at all nodes of the tree, but only on those on the upper levels. In Section \ref{subsec:5lb} we show how applying these lower bounds in some of the nodes of the search tree affects the performance of the branch-and-bound algorithm.

\newpage

\subsection{Pseudo - algorithm}\label{subsec:pseudoAlg}

We now introduce/review some notation:
\begin{center}
	\begin{tabular}{l|l}
		$N_0$ & root node	\\
		$\Ll$ & set of live nodes (maintained as a stack)\\
		${r}^*$, ${x}^*$ & radius and center of the best solution found at any stage\\		
		$\Pp(N)$ & set of points corresponding to nodes on the path from $N_0$ to $N$ (included)\\
		$\Ss t(N)$ & set of points corresponding to nodes on the subtree of $N$ (excluded)\\
		$r(N)$, $x(N)$ & radius and center of $MEB(N)$, respectively\\
		$b(N)$ & number of branches/child nodes of $N$\\
		$C(N)$ & list of child nodes of $N$
	\end{tabular}
\end{center}

We present the pseudo-code of the branch-and-bound algorithm in Algorithm \ref{alg:BB}. Calculating the lower bound on lines 7-10 is optional.

\begin{algorithm}
	\caption{A branch-and-bound algorithm for the $MkEB$ problem}\label{alg:BB}
	\begin{algorithmic}[1]
		\Require $\Pp$, and an initial solution $x^*, r^*$
		\Ensure $x^*, r^*$
		\State $\Ll =\{N_0\}.$
		
		\smallskip	
		\While{$\Ll \neq\emptyset$}
		\State Select the top node $N$ from $\Ll$.
		\If{$r(N) \geq r^*$}
			\State prune $N$; \textbf{break}
		\EndIf
		
		\MyStateA{Get $\overline{r}$, by solving the relaxation of (\ref{lbSOCPnode}) or (\ref{lbQPnode}).}
		\If{$\overline{r}\geq r^*$}
			\State prune $N$. \textbf{break}.
		\EndIf
		\MyStateA {Branch $N$:  for each $p\in \Ss t(N)$, calculate $\|p-x(N)\|$, and select the $b(N)$ points of $\Ss t(N)$ with the highest distances, sorting them in decreasing order, to get $C(N) = \{C_1, C_2 \dots, C_{b(N)}\}$, where $C_j$ corresponds to the point of $\Ss t(N)$ with the $j$-th highest distance. }
		
		\smallskip
		\ForAll {node $C_j$ in $C(N)$}
		
		\MyStateB {Get $r(C_j)$ and $x(C_j)$ by solving $MEB(\Pp(C_j))$ (see Algorithm \ref{alg:subprob}).} 
		\If {$C_j$ was not pruned} 
			\State Add $C_j$ to the top of $\Ll$.
		\EndIf
		\EndFor		
		\EndWhile
	\end{algorithmic}
\end{algorithm}

\section{Solving the MEB problem at each node}\label{sec:MEB}
Suppose node $N$ was selected to be branched and $C_1, ..., C_{b(N)}$ are its child nodes. We know the solution of the subproblem corresponding to $N$, $MEB(\Pp(N))$, and now, for each node $C_j$, we need to solve $MEB(\Pp(C_j))$. The latter subproblem differs from $MEB(\Pp(N))$ by an additional point in the set to enclose. An algorithm for the $MEB$ problem that is able to take advantage of that fact and compute with small {additional} work the solution of each successor's subproblem is crucial. A dual algorithm is consequently a natural choice, since the solution of $MEB(\Pp(N))$ is dual feasible for the problem $MEB(\Pp(C_j))$.
 
\smallskip 
 
Dearing and Zeck proposed in \cite{Dearing09} a dual algorithm based on the fact that the minimum enclosing ball of $\A$, a set of points in the Euclidean space, is determined by an affinely independent subset $\Ss$ of $\A$, called the \emph{support set}, whose points are on the boundary of the optimal ball. At each iteration of the algorithm, the current ball is optimal (namely it has minimum radius) with respect to all points it contains. However, it may not be feasible since it may not include all the points in $\A$. The algorithm is initialized with a support set $\Ss$ of two points and the trivial solution, $(x, r)$, to $MEB(\Ss)$. Then, it iterates between solving $MEB(\Ss)$ and updating $\Ss$ until feasibility is achieved. 
At each iteration, if $(x, r)$ is not optimal (feasible) for $MEB(\A)$, a point $p\in \A$, that is not yet enclosed, is chosen to enter~$\Ss$. $\Ss$ is then updated by either the addition of $p$ to $\Ss$ or by the replacement of an existing point in $\Ss$ by $p$. In either case the points in $\Ss$ remain affinely independent and the size of $\Ss$ never exceeds $n+1$. Given the updated set $\Ss$, $MEB(\Ss)$ is then solved and $p$ becomes enclosed (while some of the previously enclosed points may not be anymore). During the solution of $MEB(\Ss)$ some points may be dropped from $\Ss$. Since $MEB(\Ss)$ is a relaxation of $MEB(\A)$, the radius of $MEB(\Ss)$ is a lower bound on the optimal radius of $MEB(\A)$. The radius strictly increases at each iteration, so the algorithm terminates in a finite number of iterations.

In \cite{Cavaleiro17} the authors present a modification of Dearing and Zeck's algorithm that makes it more efficient by reducing the computational cost of each iteration from $\bigO(n^3)$ to~$\bigO(n^2)$. 

\medskip

At each node of the search tree, we use Dearing and Zeck's dual algorithm to solve the MEB problem associated to that node. To see why this dual algorithm is most suitable, suppose node~$N$ is selected to be branched. When $N$ is analyzed, $MEB(\Pp(N))$ is solved, and the corresponding radius, center, and support set information kept in the stack along with the node itself. Once $N$ is branched, problem $MEB(\Pp(C_j))$ needs to be solved for each child node $C_j$, where $\Pp(C_j)=\Pp(N) \cup \{q_j\}$ and $q_j$ is the point $C_j$ corresponds to. If $q_j$ is covered by the ball of $MEB(\Pp(N))$, then the solution of $MEB(\Pp(C_j))$ is the same as of $MEB(\Pp(N))$. Otherwise, we apply Dearing and Zeck's algorithm starting with the radius, center, and support set of $MEB(\Pp(N))$, and $q_j$ as the point to enter the support set. Typically, solving $MEB(\Pp(N) \cup \{q_j\})$ requires a modest number of iterations. Moreover, since at each iteration of the algorithm the radius increases, if it ever becomes larger than the best bound we do not need to run the algorithm until optimality; the node $q_j$ and its successors can be pruned. Finally, the improvement proposed in \cite{Cavaleiro17} has a substantial impact in the context of high dimensional MkEB instances, making the solution of $MEB(N)$ much faster for each node~$N$.

\subsection{Lower Bound on the radius of a node}\label{subsec:LBMEB}

Before solving problem $MEB(\Pp(C_j))$, we can obtain a lower bound on its optimal radius, and, if it is worse than the best solution found, use it to prune $C_j$. If $q_j$, the point corresponding to $C_j$, is not covered by $MEB(\Pp (N))$, then it is always going to be part of the optimal support set of $MEB(\Pp(N) \cup \{q_j\})$ \cite{Dearing09}, and therefore will be on the boundary of the optimal ball. Using the triangle inequality, we conclude that

\begin{equation*}
\|q_j-p\|  \leq \|p-x(C_j)\| + \|q_j-x(C_j)\| \leq 2r(C_j)\quad \forall p\in \Pp(N),
\end{equation*} 

\noindent where $x(C_j)$ and $r(C_j)$ are the center and radius of $MEB(\Pp(N) \cup \{q_j\})$, respectively. Thus, we obtain the following lower bound on $r(C_j)$

\begin{equation}\label{lb}
r(C_j) \geq  \frac 12 \max_{p\in \Pp(N)} {\left\|q_j-p\right\|}.
\end{equation} 

\noindent The burden of finding lower bound (\ref{lb}) can be avoided by calculating an easy upper bound on $r(C_j)$. This upper bound is again a consequence of the observation that, if not covered already, $q_j$ will be on the boundary of $MEB(\Pp(N) \cup \{q_j\})$. Therefore the diameter of the optimal ball will always be smaller than the distance between $q_j$ and the farthest point on the boundary of $MEB(\Pp(N))$. That distance is $\|q_j-x(N)\| + r(N)$. 
As a consequence, we obtain the following upper bound on~$r(C_j)$

\begin{equation}\label{up}
r(C_j)\leq \frac12\left(\|q_j-x(N)\| + r(N)\right).
\end{equation}

\medskip

Algorithm \ref{alg:subprob} summarizes the procedure to solve the subproblem at each node.

\begin{algorithm}
	\caption{Solving $MEB(\Pp(C_j))$ for each child node of $N$}\label{alg:subprob}
	\begin{algorithmic}[1]
		\Require $x(N), r(N), q_j, r^*$
		\Ensure $r(C_j), x(C_j)$
		\smallskip
		\State Calculate $d_j = \|q_j-x(N)\|$.
		
		\If {$d_j\leq r(N)$} 
			\State $q_j$ is already covered, so $r(C_j) = r(N)$ and $x(C_j) = x(N)$. 	
		\ElsIf {$\frac12(d_j + r(N))> {r}^*$}
			\If{$\max_{p\in \Pp(N)} \frac12\left\|q_j-p\right\| \geq r^*$}
				\State Prune $C_j$. 
			\Else
				\State Find $r(C_j)$ and $x(C_j)$ by solving $MEB(\Pp(N) \cup \{q_j\})$ using the solution of $MEB(\Pp(N))$.
			\EndIf
		\EndIf
	\end{algorithmic}
\end{algorithm}

\section{Experimental results}\label{sec:experiments}
In this section we examine the effectiveness of our algorithm by evaluating its performance on various synthetic datasets. The algorithm was implemented using MATLAB 2014a. The tests were run on a PC with an Intel Core i5 2.30GHz processor, with 4GB RAM, running Windows 7.

\smallskip

We generated our datasets with the following approaches:
\begin{enumerate}
	\item \emph{ball}: points uniformly sampled on a unit ball;
	\item \emph{ring}: points uniformly sampled on a ring with inner radius $0.8$ and outer radius $1.2$;
	\item \emph{normal}: points sampled from a multivariate standard normal distribution;
	\item \emph{exponential}: each coordinate of the points independently sampled from an exponential distribution with mean equal to 1;
	\item \emph{$b$-outliers}: points uniformly sampled on a unit ball and $b$ artificial outliers uniformly sampled from the ring with same center and radius between $1$ and $3$.
\end{enumerate}

\begin{center}
\begin{figure}[ht]
	\centering
	\subfigure[\emph{ball}]{\includegraphics[width=4.9cm]{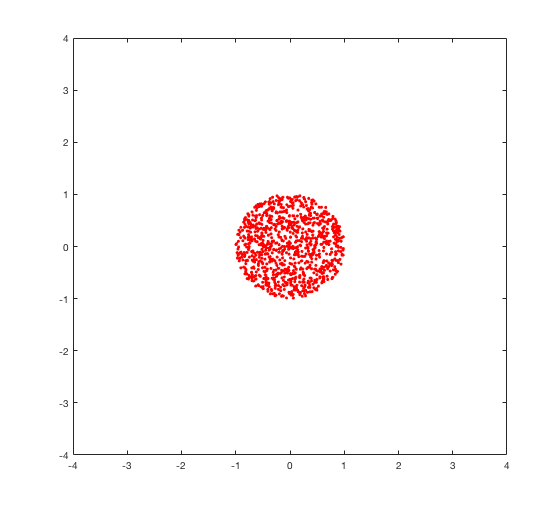}}
		\hfill
	\subfigure[\emph{ring}]{\includegraphics[width=4.9cm]{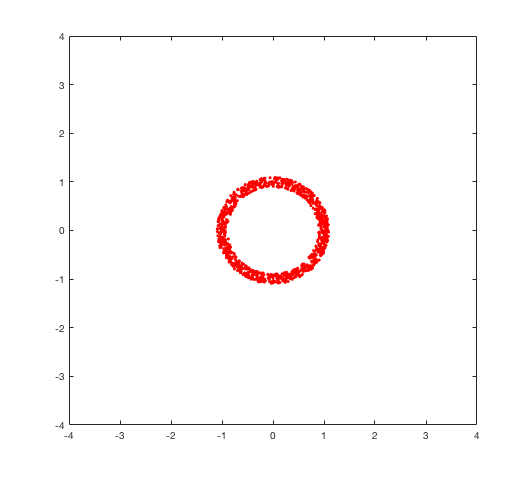}}
	\hfill
	\subfigure[\emph{normal}]{\includegraphics[width=4.9cm]{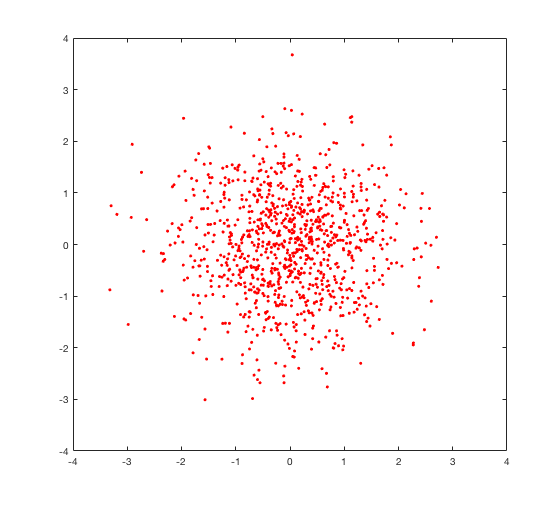}}

	\subfigure[\emph{exponential}]{\includegraphics[width=4.9cm]{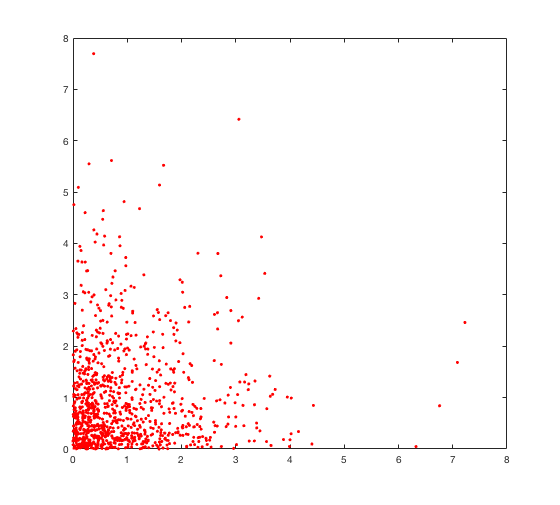}}
	\hfill
	\subfigure[\emph{$10$-outliers}]{\includegraphics[width=4.9cm]{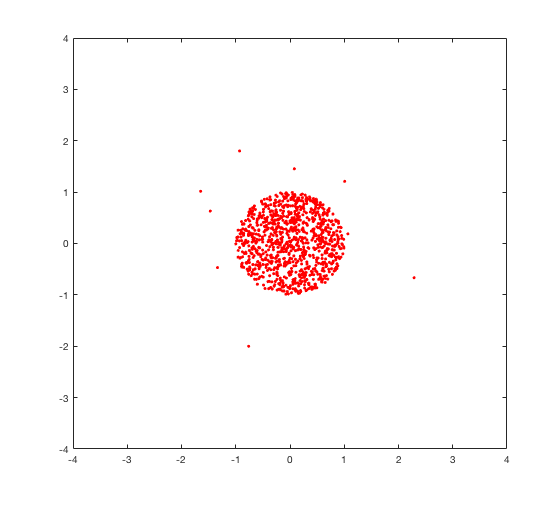}}

	\caption{Examples of the five different types of datasets with $1000$ points and dimension $2$.}\label{fig:datasets}
\end{figure}
\end{center}

Figure \ref{fig:datasets} shows $2$-dimensional instances for each of these datasets.

\subsection{Discussion on the B\&B algorithm performance}\label{subsec:5gen}

We ran the B$\&$B algorithm on different datasets, varying the dimension and the number of points. For each dataset, we generated $10$ random instances and solved the MkEB problem for different values of $k$. In these tests we did not take advantage of lower bounds to further prune the tree (so, relaxations (\ref{lbSOCPnode}) or (\ref{lbQPnode}) were not computed). For the reasons discussed in Section \ref{subsec:initSol}, to obtain an initial solution for the datasets \emph{normal} and \emph{ball}, we used the Spherical Ordering method; for the \emph{ring} dataset we selected a random point and the $k-1$ closest points to it, and found the MEB of that set; and finally, for the \emph{exponential} dataset we used Spherical Peeling. 

The following performance measures for the algorithm are reported:
\begin{itemize}
	\item \textbf{Explored nodes (EN)}: the number of explored nodes, that is, nodes where a MEB problem was solved;
	\item \textbf{\% EN optimum is found}: the percentage of explored nodes after which the optimal solution was found;
	\item \textbf{Dual iter. per node}: the ratio between the number of dual MEB iterations and number of explored nodes;
	\item \textbf{Time}: the total elapsed CPU time in seconds.	
\end{itemize}
\noindent Tables \ref{tab:m1000n2}, \ref{tab:m1000n10}, and	\ref{tab:m100n10} report the average values of the above performance measures for different experiments. 

\smallskip

The most challenging instances in terms of number of explored nodes are those that contain many $k$-subsets whose MEB radius is close to the radius of the optimal $k$-enclosing ball. In these instances, a large portion of the tree has to be scanned before finding the optimal solution and prove its optimality. The \emph{ball} instances for all values of $k$ and the \emph{ring} instances for large enough $k$ are such examples. On the other hand, when the data has many outlying points, like in the \emph{normal} or the \emph{exponential} case, the algorithm tends to explore a smaller number of nodes, because the B$\&$B will place many of the outlying points on the left side of the tree at the early stages of the algorithm. The nodes corresponding to those points are quickly realized that they can be pruned, and so are large parts of the tree. The same behavior can be observed in the \emph{ring} instances when $k$ is small, since many $k$-subsets will have points at opposite sides of the ring. 
We point out that the number of MEB iterations per explored node is close to $1$. This validates the choice of Dearing and Zeck's dual algorithm for solving MEB subproblems for each node. As we will see, while the number of MEB iterations tends to increase with the dimension, it generally remains rather small.

\smallskip

Small increases in the dimension cause a substantial increase in the difficulty of the problem due to the ``curse of dimensionality''. Consider when a node is branched and its minimum ball is inflated in order to cover a point corresponding to one of its child nodes. When the dimension is small this inflated ball has a good chance of covering many other points in the vicinity. This is true especially if the point to cover is far from the current center. When the dimension is large, however, the inflated ball has less chance to cover many other points unless the number of points also increases (exponentially) with the dimension. In other words, by increasing the dimension we have a lot more room to disperse the points and so the chance of an inflating ball to capture many other points decreases. 

Another issue with increasing the dimension is that the number of dual iterations per node increases. When a ball is inflated and translated to enclose a new point, as the dimension increases, the chance that a previously covered point becomes uncovered also increases. For the same reason, this behavior tends to become more frequent for larger $k$. 

Table~\ref{tab:m1000n10} reports the performance results for the \emph{normal} case, with $n=10$ and $k=50$, $k=900$, and $k=990$. These can be compared with the analogous results reported in Table~\ref{tab:m1000n2} to show how the dimension increases the complexity of the problem. The number of explored nodes in dimension $10$ increases significantly with respect to dimension $2$, and for larger values of $k$, the average number of dual iterations per node also increases modestly.

\smallskip

\begin{figure}[htbp]
	\centering
	
	\subfigure[Average (dual) iterations as a function of~$k$]{\includegraphics[width=7cm]{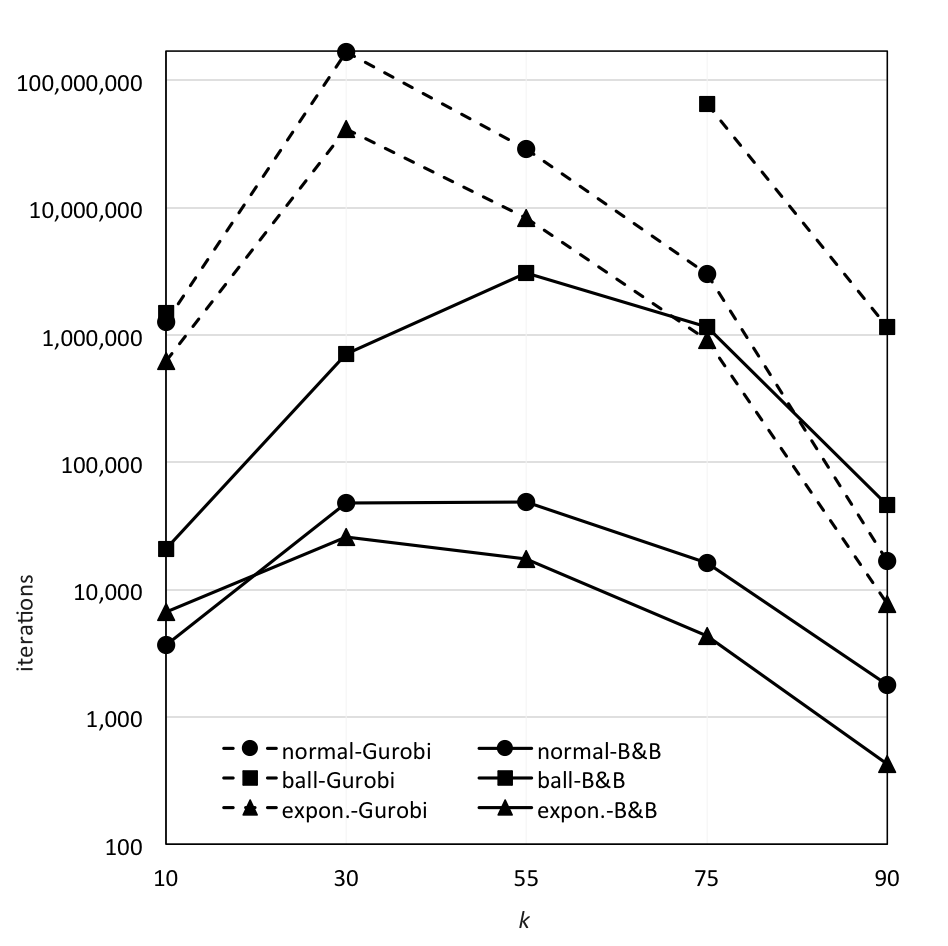}}
	\hfill
	\subfigure[Average CPU time in seconds as a function of~$k$]{\includegraphics[width=7cm]{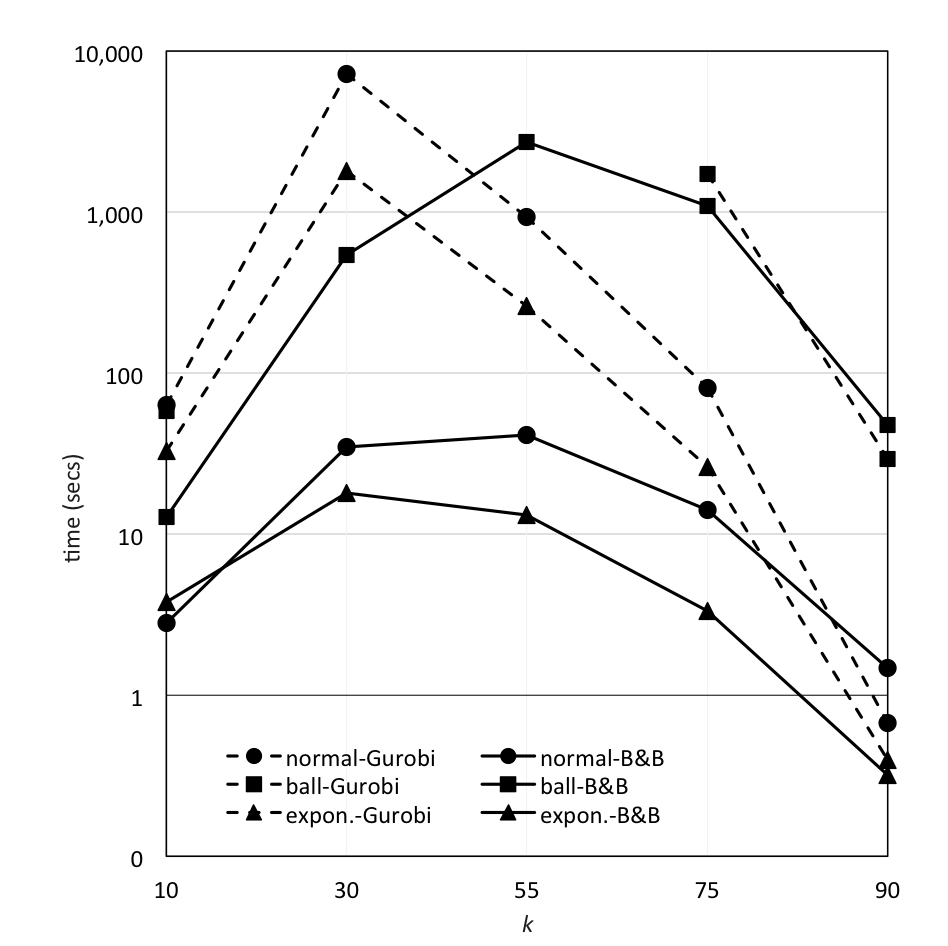}}
	\hfill
	\caption{Gurobi vs. B$\&$B for different $10$-dimensional \emph{normal}, \emph{ball}, and \emph{exponential} datasets with $100$ points. The data for the \emph{ball} dataset for $k=30$ and $k=55$ is missing since Gurobi was not able to solve these problems to optimality in the allowed time.}
	\label{fig:gurobi}
\end{figure}

We compared the results of Table~\ref{tab:m100n10} to those of the general-purpose MIQP-solver of Gurobi \cite{gurobi} (version 6.5.2) using its MATLAB interface. All Gurobi parameters were kept at their default values\footnote{We realize that by fine-tuning these parameters Gurobi may have a much better performance, but the idea here is to show that our special purpose algorithm works better than a generic and default branch-and-bound method using the dual QP Simplex algorithm.}, with only the following exception: the dual simplex method was the algorithm chosen for solving the relaxation at each node. Moreover, the same initial solutions given in the experiments of Table~\ref{tab:m100n10} were given to Gurobi as warm start. For each dataset we compared the number of Gurobi's dual simplex iterations and the number of Dearing and Zeck's dual algorithm iterations in the B$\&$B (Fig. \ref{fig:gurobi} (a)); we also compared the CPU time of Gurobi and the B$\&$B  (Fig. \ref{fig:gurobi} (b)). The maximum time allowed was $10^5$ seconds for both algorithms.

We observe that our method was faster for all the datasets studied. It is of course not surprising that a dedicated algorithm is superior to a general-purpose code. Still, the comparison is necessary in order to argue that off-the-shelf methods cannot successfully compete with our approach. Our comparison is based on the $100$-point instances, since Gurobi was not able to solve the instances with $1000$ points within the allocated time.

\smallskip

In order to show that larger instances can also be solved in useful time, for certain datasets, we run some experiments for $10000$-point instances of \emph{$b$-outliers} datasets in dimension~$10$ and $100$, for different values of $b$. Table~\ref{tab:outliers} shows the average results of $10$ runs.

\begin{table}[htbp]	
	\footnotesize \centering\setlength{\tabcolsep}{6pt}
	\begin{tabular}{lr r rrrr}
		\noalign{\smallskip}\hline \noalign{\smallskip}
		 &  && \multicolumn{1}{c}{Explored} & \multicolumn{1}{c}{\% EN optimum} & \multicolumn{1}{c}{Dual iter.}  &\\
		Dataset & \multicolumn{1}{c}{$k$} && \multicolumn{1}{c}{nodes (EN)} & \multicolumn{1}{c}{is found} & \multicolumn{1}{c}{per node}  & \multicolumn{1}{c}{Time} \\
		\noalign{\smallskip}\cline{1-2}\cline{4-7}\noalign{\smallskip}
		
		\emph{normal} & 50 && 231243 &  11\%  & 1.00  & 150.12 \\
		\emph{normal} & 100 && 315868 &  19\% & 1.00  & 220.39 \\
		\emph{normal} & 250 && 717512 &  42\%  & 1.01  & 580.14 \\
		\emph{normal} & 750 && 228192 &  40\%  & 1.01  & 232.68 \\
		\emph{normal} & 900 && 24439 &  34\% & 1.02  & 30.60 \\
		\emph{normal} & 990 && 94    &  44\% & 1.05  & 0.15 \\
		\noalign{\smallskip}
		\emph{ball} & 50 && 307411 &  43\% & 1.00  & 228.58 \\
		\emph{ball} & 100 && 797306 &  69\%  & 1.01  & 638.79 \\
		\emph{ball} & 250 && 8295548 &  74\%  & 1.00  & 7401.77 \\
		\emph{ball} & 750 && 3933869 &  72\%& 1.00  & 4875.08 \\
		\emph{ball} & 900 && 543560 &  57\% & 1.01  & 518.71 \\
		\emph{ball} & 990 && 772   &  65\% & 1.11  & 0.87 \\
		\noalign{\smallskip}
		\emph{ring} & 50 && 185739 &  91\%& 1.00  & 160.55 \\
		\emph{ring} & 100 && 347188 &  80\%  & 1.00  & 295.75 \\
		\emph{ring} & 250 && 1239892 &  87\%  & 1.00  & 1152.42 \\
		\emph{ring} & 750 && 25015147 &  64\% & 1.00  & 19527.62 \\
		\emph{ring} & 900 && 1591675 &  70\% & 1.01  & 1204.24 \\
		\emph{ring} & 990 && 1082  &  46\%& 1.12  & 0.90 \\
		\noalign{\smallskip}
		\emph{exponential} & 50 && 275855 &  98\% & 1.00  & 188.83 \\
		\emph{exponential} & 100 && 400178 &  99\% & 1.00  & 267.75 \\
		\emph{exponential} & 250 && 532264 &  94\% & 1.00  & 356.61 \\
		\emph{exponential} & 750 && 568182 &  95\% & 1.00  & 439.69 \\
		\emph{exponential} & 900 && 47783 &  97\% & 1.02  & 4 16 \\
		\emph{exponential} & 990 && 130   &  93\%& 1.16  & 0.16 \\
		\noalign{\smallskip}\hline
	\end{tabular}	
	\caption{Performance of the B$\&$B for different $2$-dimensional datasets with $1000$ points}\label{tab:m1000n2}		
\end{table}

\begin{table}[htbp]
	\footnotesize \centering\setlength{\tabcolsep}{6pt}
	\begin{tabular}{lr r rrrr}
		\noalign{\smallskip}\hline \noalign{\smallskip}
		 &  && \multicolumn{1}{c}{Explored} & \multicolumn{1}{c}{\% EN optimum} & \multicolumn{1}{c}{Dual iter.}  &\\
		 Dataset & \multicolumn{1}{c}{$k$} && \multicolumn{1}{c}{nodes (EN)} & \multicolumn{1}{c}{is found} & \multicolumn{1}{c}{per node}  & \multicolumn{1}{c}{Time} \\
		\noalign{\smallskip}\cline{1-2}\cline{4-7}\noalign{\smallskip}
		
		\emph{normal} & 50    && 33906317& 37\%  & 1.03  & 36503.59 \\
		\emph{normal} & 900   && 450714331 & 63\%  & 1.06  & 575461.26 \\
		\emph{normal} & 990   && 3929 & 55\%  & 1.36  & 9.89 \\
		\noalign{\smallskip}\hline
	\end{tabular}	
	\caption{Performance of the B$\&$B for $10$-dimensional  \emph{normal} datasets with $1000$ points.}\label{tab:m1000n10}
\end{table}

\begin{table}[htbp]
	\footnotesize \centering\setlength{\tabcolsep}{6pt}
	\begin{tabular}{lr r rrrr}
		\noalign{\smallskip}\hline \noalign{\smallskip}
		 &  && \multicolumn{1}{c}{Explored} & \multicolumn{1}{c}{\% EN optimum} & \multicolumn{1}{c}{Dual iter.}  &\\
		 Dataset & \multicolumn{1}{c}{$k$} && \multicolumn{1}{c}{nodes (EN)} & \multicolumn{1}{c}{is found} & \multicolumn{1}{c}{per node}  & \multicolumn{1}{c}{Time} \\
		\noalign{\smallskip}\cline{1-2}\cline{4-7}\noalign{\smallskip}
		
		\emph{normal} & 10    && 3426  & 38\%  & 1.06  & 2.81 \\
		\emph{normal} & 30    && 42847 & 60\%  & 1.11  & 34.66 \\
		\emph{normal} & 55    && 43179 & 43\%  & 1.13  & 41.43 \\
		\emph{normal} & 75    && 13519 & 48\%  & 1.19  & 14.22 \\
		\emph{normal} & 90    && 1315  & 57\%  & 1.35  & 1.47 \\
		\noalign{\smallskip}
		\emph{ball} & 10    && 19973 & 66\%  & 1.05  & 12.76 \\
		\emph{ball} & 30    && 662078 & 71\%  & 1.08  & 544.00 \\
		\emph{ball} & 55    && 2775139 & 88\%  & 1.11  & 2713.10 \\
		\emph{ball} & 75    && 973051 & 60\%  & 1.19  & 1094.62 \\
		\emph{ball} & 90    && 32722 & 63\%  & 1.41  & 47.97 \\
		\noalign{\smallskip}
		\emph{exponential} & 10  && 6396  & 66\%  & 1.04  & 3.80 \\
		\emph{exponential} & 30    && 23390 & 79\%  & 1.10  & 18.04 \\
		\emph{exponential} & 55    && 15183 & 80\%  & 1.16  & 13.09 \\
		\emph{exponential} & 75    && 3583  & 69\%  & 1.20  & 3.35 \\
		\emph{exponential} & 90    && 298   & 53\%  & 1.42  & 0.32 \\
		\noalign{\smallskip}\hline
	\end{tabular}	
	\caption{Performance of the B$\&$B for different $10$-dimensional datasets with $100$ points.}\label{tab:m100n10}
\end{table}

\begin{table}[htbp]	
	\footnotesize \centering\setlength{\tabcolsep}{6pt}
	\begin{tabular}{lclr r rrrr}
		\noalign{\smallskip}\hline \noalign{\smallskip}
		& & &&& \multicolumn{1}{c}{Explored} & \multicolumn{1}{c}{\% EN optimum} & \multicolumn{1}{c}{Dual iter.}  &\\
		 Dataset & \multicolumn{1}{c}{$n$} &\multicolumn{1}{c}{$b$} & \multicolumn{1}{c}{$k$} && \multicolumn{1}{c}{nodes (EN)} & \multicolumn{1}{c}{is found} & \multicolumn{1}{c}{per node}  & \multicolumn{1}{c}{Time} \\
				 
		\noalign{\smallskip}\cline{1-4}\cline{6-9}\noalign{\smallskip}
		\emph{$b$-outliers} & 10    & 10 (0.1\%) & 9990  && 69    & 84\%  & 3.46  & 0.72 \\
		\emph{$b$-outliers} & 10    & 50 (0.5\%) & 9950  && 3593  & 96\%  & 2.74  & 37.40 \\
		\emph{$b$-outliers} & 10    & 100 (1\%) & 9900  && 43012 & 98\%  & 2.43  & 459.19 \\
		\noalign{\smallskip}
		\emph{$b$-outliers} & 100   & 10 (0.1\%) & 9990  && 58    & 57\%  & 15.29 & 17.22 \\
		\emph{$b$-outliers} & 100   & 50 (0.5\%) & 9950  && 814   & 52\%  & 3.03  & 46.62 \\
		\emph{$b$-outliers} & 100   & 100 (1\%) & 9990  && 9711  & 55\%  & 1.71  & 245.22 \\
		\noalign{\smallskip}\hline
	\end{tabular}
	\caption{Performance of the B$\&$B for $10000$-point \emph{$b$-outliers} datasets for different values of $b$ and dimension $n$.}\label{tab:outliers}
\end{table}

\newpage

\subsection{On the impact of using lower bounds}\label{subsec:5lb}

In Section \ref{subsec:LB} we introduced two relaxations of the problem of finding the best $k$-enclosing ball on the subtree of a node: (\ref{lbSOCPnode}) is an SOCP and (\ref{lbQPnode}) is a QP. These methods can be applied on some or all nodes of the search tree in order to prune parts of the tree more aggressively. 

We carried out some experiments to understand how these relaxations help us reduce the number of explored nodes. Both relaxations were applied separately for each instance. We did not compute these lower bounds on every single node, but only when all the following were met: the level of the node was greater than $1$; the number of points involved in the relaxation was larger than $k$ (otherwise the relaxation is simply an MEB problem); and the number of points on the subtree of a node (that is the set $\Ss t(N)$) was at least $0.1N$. Note than when the number of points in $\Ss t(N)$ is small, it is likely that the radius returned by the relaxation is not going to be much larger than that of $MEB(\Pp (N))$. Of course, the relaxation was never solved on leaf nodes. 

\smallskip

The experiments were conducted under the same conditions as the experiments of Tables~\ref{tab:m1000n2} and~\ref{tab:m100n10}. The results are presented in Tables \ref{tab:LB1000} and \ref{tab:LB100}, with the following measures
\begin{itemize}
	\item \textbf{EN {SOCP} (\%)}: percent ratio of the number of explored nodes when the SOCP relaxation (\ref{lbSOCPnode}) was used, over the number of explored nodes when no lower bound was used;
	\item \textbf{EN {QP} (\%)}:  percent ratio of the number of explored nodes when the QP relaxation (\ref{lbQPnode}) was used, over the number of explored nodes when no lower bound was used.
\end{itemize}

We observe that the QP relaxation had almost no impact on decreasing the number of explored nodes. This phenomenon has to do with the fact that the value of $M$ is usually very large. This, in turn, affects the optimal relaxed solution and causes it to be not much larger than the radius of $MEB(\Pp (N))$. A small reduction was however observed when $k$ is very large, which is explained by the fact that, since the value of most of the $\beta_j$ will be close to 1, the effect of $M$ is neutralized for most of the constraints. The SOCP relaxation performs better. However, as mentioned earlier, this relaxation is full-fledged SOCP and at this point no effective algorithm other than interior point methods are known to solve it. As a result, it becomes exceedingly time-consuming to apply it to every node. 

In general terms, it is not surprising that these relaxations have less impact when $k$ is small. Note that the constraint $\sum \beta_j \geq k-l(N)$, when $k$ is small, causes the solution of relaxations to be too small (compared to the value of $MEB(\Pp (N))$) to produce useful values for pruning. 

As a final comment, the (usually low) impact of using these relaxations cannot be only  justified by formulation features. We must also consider the fact that due to the node/point attribution rule at each node, the radius of the MEB of $\Pp(N)$ is already a very good lower bound and blunts the effect of these relaxations.  In fact, these relaxations, when applied to a branch-and-bound based on the lexicographic node/point attribution rule, performed much better and resulted in much more aggressive pruning.

\begin{table}[htbp]	
	\footnotesize \centering\setlength{\tabcolsep}{6pt}
	\begin{tabular}{lr r rr}
		\noalign{\smallskip}\hline \noalign{\smallskip}
		Dataset & \multicolumn{1}{c}{k} && EN {SOCP} (\%) & EN {QP} (\%) \\
		\noalign{\smallskip}\cline{1-2}\cline{4-5}\noalign{\smallskip}			
		\emph{normal} & 50    && 96\%  & 100\% \\
		\emph{normal} & 100   && 85\%  & 100\% \\
		\emph{normal} & 250   && 84\%  & 99\% \\
		\emph{normal} & 750   && 89\%  & 100\% \\
		\emph{normal} & 900   && 85\%  & 99\% \\
		\emph{normal} & 990   && 91\%  & 96\% \\
		\noalign{\smallskip}
		\emph{exponential} & 50    && 80\%  & 99\% \\
		\emph{exponential} & 100   && 61\%  & 97\% \\
		\emph{exponential} & 250   && 49\%  & 98\% \\
		\emph{exponential} & 750   && 74\%  & 99\% \\
		\emph{exponential} & 900   && 86\%  & 99\% \\
		\emph{exponential} & 990   && 93\%  & 95\% \\
		\noalign{\smallskip}\hline \noalign{\smallskip}
	\end{tabular}
	\caption{Effect of relaxations (\ref{lbSOCPnode}) and (\ref{lbQPnode}) on the number of explored nodes, for $2$-dimensional datasets with $1000$ points.}\label{tab:LB1000}
\end{table}

\begin{table}[htbp]	
	\footnotesize \centering\setlength{\tabcolsep}{6pt}
	\begin{tabular}{lr r rr}
		\noalign{\smallskip}\hline \noalign{\smallskip}
		Dataset & \multicolumn{1}{c}{k} &&  EN {SOCP} (\%) & EN {QP} (\%) \\
		\noalign{\smallskip}\cline{1-2}\cline{4-5}\noalign{\smallskip}		
		\emph{normal} & 10    &       & 99\%  & 99\% \\
		\emph{normal} & 30    &       & 91\%  & 99\% \\
		\emph{normal} & 55    &       & 89\%  & 99\% \\
		\emph{normal} & 75    &       & 90\%  & 99\% \\
		\emph{normal} & 90    &       & 93\%  & 99\% \\
		\noalign{\smallskip}
		\emph{exponential} & 10    &       & 95\%  & 97\% \\
		\emph{exponential} & 30    &       & 87\%  & 98\% \\
		\emph{exponential} & 55    &       & 85\%  & 98\% \\
		\emph{exponential} & 75    &       & 86\%  & 98\% \\
		\emph{exponential} & 90    &       & 89\%  & 96\% \\
		\noalign{\smallskip}\hline \noalign{\smallskip}
	\end{tabular}
	\caption{Effect of relaxations (\ref{lbSOCPnode}) and (\ref{lbQPnode}) on the number of explored nodes, for $10$-dimensional datasets with $100$ points.}	\label{tab:LB100}
\end{table}

\subsection{On the impact of an initial solution}\label{subsec:5is}

Starting the B\&B algorithm with the knowledge of a good upper bound on the optimal solution can be potentially advantageous. In Section \ref{subsec:initSol} we presented different methods to obtain such initial solutions. But how much effort should be put in finding a good initial solution? We performed some experiments  to see how much of an impact good initial solutions can have on the performance of the algorithm. In particular, we wished to see whether or not a good initial solution is essential. We conducted two diametrically opposite experiments: In one we ran our B$\&$B algorithm without any initial solution. In the second we used the optimal solution (obtained in an earlier run) as the initial solution. Both methods were tested on the same datasets. We used the following performance measures:

\begin{itemize}
	\item \textbf{EN optimum given as initial sol. (\%)}: percent ratio of the number of explored nodes when the optimum was given as initial solution, over the number of explored nodes, when no initial solution was provided;
	\item \textbf{EN optimum found with no initial sol. (\%)}: the percentage of the number of explored nodes after which the optimal solution was found, over the total, when no initial solution was provided to the algorithm. 
\end{itemize}

\noindent Tables~\ref{tab:initSol1000} and~\ref{tab:initSol100} contain the averages of the measures above for different 10 samples of each dataset.

\smallskip

Two factors explain why, in some cases, the number of explored nodes when the optimum is provided as initial solution is close to the number of explored nodes when no initial solution is given: (a) either a good upper bound or the optimum is found early, when no initial solution is given, or (b) after finding a good upper bound, or even the optimal solution, the algorithm still needs to scan a large part of the search tree to prove it has found the optimum. It is also possible that both these factors are at work in some cases. For example, in the \emph{normal} datasets, when no initial solution is given, the optimal solution is found earlier for smaller $k$. However, the reduction in the number of explored nodes is comparable for different values of $k$. So factor (a) seems to have more weight for smaller $k$ while factor (b) has more weight for larger $k$. 

On the other hand, if the number of explored nodes when the optimum is provided as initial solution is significantly smaller than when no initial solution is given, then, when no initial solution is given, a good upper bound is found in a later stage. In this case the knowledge of a good initial solution makes a considerable difference. This can be observed in the \emph{exponential} dataset with $n=2$ and $m=1000$.

The main conclusion, therefore, is that although an initial solution can sometimes reduce the number of explored nodes considerably, in many cases it does not. Therefore, spending too much effort at the outset of the algorithm to find a good upper bound can often be wasteful.

\begin{table}[htbp]	
	\footnotesize \centering\setlength{\tabcolsep}{6pt}
	\begin{tabular}{lr r rr}
		\noalign{\smallskip}\hline \noalign{\smallskip}
		& &       												& EN optimum given & EN optimum found \\
		Dataset & \multicolumn{1}{c}{k} &       & as initial sol. (\%)  & with no initial sol. (\%)\\
		\noalign{\smallskip}\cline{1-2}\cline{4-5}\noalign{\smallskip}		
		\emph{normal} & 50    &       & 96\%  & 14\% \\
		\emph{normal} & 100   &       & 94\%  & 22\% \\
		\emph{normal} & 250   &       & 90\%  & 37\% \\
		\emph{normal} & 750   &       & 94\%  & 40\% \\
		\emph{normal} & 900   &       & 94\%  & 44\% \\
		\emph{normal} & 990   &       & 93\%  & 50\% \\
		\noalign{\smallskip}
		\emph{ball} & 50    &       & 90\%  & 64\% \\
		\emph{ball} & 100   &       & 81\%  & 72\% \\
		\emph{ball} & 250   &       & 86\%  & 61\% \\
		\emph{ball} & 750   &       & 94\%  & 72\% \\
		\emph{ball} & 900   &       & 96\%  & 59\% \\
		\emph{ball} & 990   &       & 94\%  & 46\% \\
		\noalign{\smallskip}
		\emph{exponential} & 50    &       & 64\%  & 97\% \\
		\emph{exponential} & 100   &       & 50\%  & 96\% \\
		\emph{exponential} & 250   &       & 71\%  & 95\% \\
		\emph{exponential} & 750   &       & 35\%  & 94\% \\
		\emph{exponential} & 900   &       & 34\%  & 98\% \\
		\emph{exponential} & 990   &       & 66\%  & 93\% \\
		\noalign{\smallskip}\hline \noalign{\smallskip}
\end{tabular}
\caption{Maximum impact of an initial solution, for $2$-dimensional datasets with $1000$ points.}\label{tab:initSol1000}
\end{table}

\begin{table}[htbp]	
	\footnotesize \centering\setlength{\tabcolsep}{6pt}
	\begin{tabular}{lr r rr}
		\noalign{\smallskip}\hline \noalign{\smallskip}
		& &       												& EN optimum given & EN optimum found \\
		Dataset & \multicolumn{1}{c}{k} &       & as initial sol. (\%)  & with no initial sol. (\%)\\
		\noalign{\smallskip}\cline{1-2}\cline{4-5}\noalign{\smallskip}		
		\emph{normal} & 10    &       & 86\%  & 33\% \\
		\emph{normal} & 30    &       & 88\%  & 37\% \\
		\emph{normal} & 55    &       & 85\%  & 46\% \\
		\emph{normal} & 75    &       & 83\%  & 51\% \\
		\emph{normal} & 90    &       & 86\%  & 52\% \\
		\noalign{\smallskip}		
		\emph{ball} & 10    &       & 69\%  & 82\% \\
		\emph{ball} & 30    &       & 67\%  & 75\% \\
		\emph{ball} & 55    &       & 68\%  & 76\% \\
		\emph{ball} & 75    &       & 76\%  & 92\% \\
		\emph{ball} & 90    &       & 78\%  & 84\% \\
		\noalign{\smallskip}		
		\emph{exponential} & 10    &       & 68\%  & 67\% \\
		\emph{exponential} & 30    &       & 58\%  & 82\% \\
		\emph{exponential} & 55    &       & 64\%  & 76\% \\
		\emph{exponential} & 75    &       & 77\%  & 61\% \\
		\emph{exponential} & 90    &       & 89\%  & 44\% \\
		\noalign{\smallskip}\hline \noalign{\smallskip}
		\end{tabular}
	\caption{Maximum impact of an initial solution, for $10$-dimensional datasets with $100$ points.}\label{tab:initSol100}
\end{table}

\section{Conclusions}\label{sec:conclusions}
In this paper we propose a branch-and-bound algorithm to solve the NP-hard problem of finding the ball of minimum radius that encloses at least $k$ points from a given set. Our algorithm makes a correspondence between the subsets of points and the tree nodes. It uses a heuristic that aims at placing subsets of points enclosed in larger balls at the root of large subtrees. As a result, it has the ability to prune very large subsets of points from the search tree early on.  This heuristic, combined with a Last-In-First-Out strategy for selecting the next node to explore, results in a combination of depth-first-search and best-first-search techniques. The algorithm retains the advantages of both these methods without the typical disadvantages. In fact, an important consequence of our LIFO approach is that the memory required to hold the active set of nodes is bounded by $m-k$.

In general, the algorithm does quite well on datasets with the presence of outliers, and the most difficult datasets for the algorithm seem to be those with an uniform distribution of the points drawn from a ball. Our experimental results also show that in almost all cases a small fraction of the nodes of the search tree need to be explored. They also support our choice of the Dearing and Zeck dual algorithm to solve the subproblems at each node, since in general we need a very small number of iterations per node. The fact that each iteration has a computational work of $\bigO(n^2)$ further speeds up the search process. 

We also studied two relaxations of the MkEB problem in order to find better lower bounds. One was based on a SOCP relaxation while another was base on a QP. Our study revealed that these lower bounds are not very effective given our node/point attribution rule, but they can be very useful if the sequences of points in the tree are to be organized in arbitrary fashion (for example in lexicographic order). 

Finally, we observed that the use of an initial good solution is not very consequential in reducing the number of explored nodes.

\bibliography{references}
\bibliographystyle{siam}

\end{document}